\documentclass[a4paper,10pt]{article}

\usepackage[ascii]{inputenc}
\usepackage{amsmath,amsfonts,amssymb,amsthm}

\newtheoremstyle{theorem}{5pt}{5pt}{\slshape}{}{\bfseries}{.}{.5em}{}

\theoremstyle{theorem}
\newtheorem{theorem}{Theorem}
\newtheorem{lemma}{Lemma}

\usepackage{titlesec}
\titleformat{\section}
{\normalfont\bfseries}{\thesection}{1em}{\large}
\titlespacing*{\section}{0pt}{3.5ex plus 1ex minus .2ex}{2.3ex plus .2ex}

\titleformat{\subsection}
{\normalfont\bfseries}{\thesubsection}{1em}{\normalsize}
\titlespacing*{\section}{0pt}{3.5ex plus 1ex minus .2ex}{2.3ex plus .2ex}

\allowdisplaybreaks[3]

\newtheorem{hypothesis}{Hypothesis}

\begin{document}

\title{Transmission Eigenvalues for a Class of Non-Compactly Supported Potentials}
\author{Esa V\!. Vesalainen}
\date{}

\maketitle

\begin{abstract}
Let $\Omega\subseteq\mathbb R^n$ be a non-empty open set for which the Sobolev embedding $H_0^2(\Omega)\longrightarrow L^2(\Omega)$ is compact, and let $V\in L^\infty(\Omega)$ be a potential taking only positive real values and satisfying the asymptotics $V(\cdot)\asymp\left\langle\cdot\right\rangle^{-\alpha}$ for some $\alpha\in\left]3,\infty\right[$. We establish the discreteness of the set of real transmission eigenvalues for both Schr\"odinger and Helmholtz scattering with these potentials.
\end{abstract}

\section{Introduction}

\subsection{Non-scattering energies and non-scattering wavenumbers}

We shall be concerned with the interior transmission problem for the Schr\"odinger and Helmholtz equations. Inverse scattering theory, and the study of the linear sampling method and the factorization method in particular, gives rise to the study of non-scattering energies. These are energies $\lambda\in\mathbb R_+$ for which there exists a non-zero incoming wave which does not scatter in the sense that the corresponding scattered wave has a vanishing main term in its asymptotic expansion. In the case of the Schr\"odinger equation with a short-range potential $V\in L^2_{\mathrm{loc}}(\mathbb R^n)$ this ultimately means that the system
\[\left\{\!\!\begin{array}{l}
\left(-\Delta+V-\lambda\right)v=0,\\
\left(-\Delta-\lambda\right)w=0
\end{array}\right.\]
has a solution $v,w\in B^\ast_2\setminus 0$ where the two functions are connected by the asymptotic condition $v-w\in\mathring{B}^\ast_2$.
The Helmholtz case is otherwise the same, except that the perturbed equation for $v$ is
\[\left(-\Delta+\lambda V-\lambda\right)\,v=0,\]
and the term non-scattering wavenumber is more appropriate.

Here the solutions are taken from the function spaces
\[B^\ast_2=\left\{u\in B^\ast\,\middle|\,\vphantom{\Big|}\partial^\gamma u\in B^\ast,\forall\left|\gamma\right|\leqslant2\right\},\]
and
\[\mathring B^\ast_2=\left\{u\in \mathring B^\ast\,\middle|\,\vphantom{\Big|}\partial^\gamma u\in \mathring B^\ast,\forall\left|\gamma\right|\leqslant2\right\},\]
where $B^\ast$ consists of those functions $u\in L_{\mathrm{loc}}^2(\mathbb R^n)$ for which
\[\sup_{R>1}\frac1R\int\limits_{B(0,R)}\left|u\right|^2<\infty,\]
where $B(0,R)$ is the ball of radius $R$ centered at the origin, and $\mathring B^\ast$ consists of those functions $u\in B^\ast$ for which
\[\frac1R\int\limits_{B(0,R)}\left|u\right|^2\longrightarrow0\]
as $R\longrightarrow\infty$. A function $V\in L^2_{\mathrm{loc}}(\mathbb R^n)$ is a short-range potential for instance when $V(\cdot)\ll\left\langle\cdot\right\rangle^{-\alpha}$ in $\mathbb R^n$ for some $\alpha\in\left]1,\infty\right[$.
For a presentation of short-range scattering theory, see e.g. Chapter XIV of \cite{HormanderII} and the first sections of \cite{PaivarintaEtSylvester}.

\subsection{Interior transmission eigenvalues}

If the potential $V$ vanishes outside a suitable bounded domain $\Omega$, then the functions $v$ and $w$ coincide outside $\Omega$ (by Rellich's lemma and unique continuation) and we are left with a solution to the problem
\[\left\{\!\!\begin{array}{l}
\left(-\Delta+V-\lambda\right)v=0,\\
\left(-\Delta-\lambda\right)w=0,
\end{array}\right.\]
where $v$ and $w$ are to be in $H_{\mathrm{loc}}^2(\Omega)\cap L^2(\Omega)$ and to satisfy the boundary conditions $v-w\in H_0^2(\Omega)$. This problem is the interior transmission problem for $\Omega$ and $V$.

Typical first steps in the study of the interior transmission problem are the finite multiplicity of transmission eigenvalues, the discreteness of the set of transmission eigenvalues, and the existence of infinitely many transmission eigenvalues.

\subsection{The purpose and the motivation of this work}

Since scattering theory does not really care about the support of $V$, it is natural to ask whether the study of the interior transmission problem can be carried over to non-compact supports. A particularly strong motivation for studying this is that, metaphorically speaking, non-scattering energies are transmission eigenvalues for the domain $\Omega=\mathbb R^n$.
In this particular case, the combinations of the techniques of short-range scattering theory and interior transmission eigenvalue problems might allow a new approach to directly deal with non-scattering energies.

One particular question which might be approached in this way is the existence of non-scattering energies. For compactly supported radial scatterers, there are always infinitely many of them as in that case the non-scattering energies coincide with the transmission eigenvalues. On the other hand, it was recently shown by Bl{\aa}sten, P\"aiv\"arinta and Sylvester \cite{BlastenEtPaivarintaEtSylvester} that for a large class of potentials there are no non-scattering energies. It is not yet known if non-scattering energies can exist for non-radial potentials.

In the following we shall take first steps into the direction of non-compact supports by considering interior transmission eigenvalues for non-compact $\Omega$ which are nearly compact in the sense that they have a suitable compact Sobolev embedding, and for potentials $V$ taking only positive real values and having a certain kind of asymptotic behaviour. For these potentials, we shall prove the basic discreteness result. This is done by proving the basic discreteness and existence results for a closely connected fourth-order equation. The more usual case of bounded $\Omega$ with a positive real-valued potential, which is bounded and bounded away from zero, is covered as a special case, including the corresponding existence result for Helmholtz transmission eigenvalues.

It should be noted that this discreteness result would imply the discreteness for the corresponding non-scattering energies if a conclusion analogous to that of Rellich's lemma could be somehow obtained. It seems that there are no known generalizations of Rellich's lemma to non-compact domains, but such generalizations might exist. We intend to return to this topic in the future.

\subsection{A few words on the preceding work}

The interior transmission problem first appeared in the papers of Kirsch \cite{Kirsch}, and Colton and Monk \cite{ColtonEtMonk}. The first papers considered radial potentials and discreteness for general potentials, see e.g. the survey \cite{ColtonEtPaivarintaEtSylvester} of Colton, P\"aiv\"arinta and Sylvester. The first general existence result was obtained by P\"aiv\"arinta and Sylvester \cite{PaivarintaEtSylvester}, and later Cakoni, Gintides and Haddar \cite{CakoniEtGintidesEtHaddar} proved the first general result on existence of infinitely many transmission eigenvalues.

It should be noted that the methods in the papers of Sylvester \cite{Sylvester}, Lakshtanov and Vainberg \cite{LakshtanovEtVainberg} and Robbiano \cite{Robbiano} are able to handle compactly supported potentials with fairly arbitrary behaviour inside the domain. I.e. the main assumptions only deal with the behaviour of the potentials in a neighbourhood of the boundary.

It is clear that we can not give here an exhaustive list of previous results and references. For a recent survey on the topic, we recommend the article \cite{CakoniEtHaddar} by Cakoni and Haddar.

The main result of this paper and its proof are in their spirit closest to the work of Hickmann \cite{Hickmann}, Serov and Sylvester \cite{SerovEtSylvester}, and Serov \cite{Serov}, who proved discreteness and existence results in compact domains for potentials exhibiting well controlled degenerate or singular behaviour at the boundary of the domain using quadratic forms, suitable weighted spaces and Hardy-type inequalities.

\subsection{On notation}

We shall employ the standard asymptotic notation. Given two complex functions $A$ and $B$ defined on some set $\Omega$, the relation $A\ll B$ means that $\left|A\right|\leqslant C\left|B\right|$ in $\Omega$ for some positive real constant $C$. The relation $A\asymp B$ means that both $A\ll B$ and $A\gg B$, and $A\gg B$ means the same as $B\ll A$. We do not insist on the implicit constants being computable.

When the letter $\varepsilon$ appears in various exponents, it denotes an arbitrarily small, and also sufficiently small, positive real constant, which usually changes its value from one occurrence to the next. The usage of this notational device should be rather transparent.

For a vector $\xi\in\mathbb R^n$, we let $\left\langle\xi\right\rangle$ denote $\sqrt{1+\left|\xi\right|^2}$, as usual.

\section{The main theorems}

We fix the dimension $n\in\mathbb Z_+$ of the ambient Euclidean space for the entire text, and all implicit constants are allowed to depend on it. Let $\Omega\subseteq\mathbb R^n$ be an open set for which the Sobolev embedding $H_0^2(\Omega)\longrightarrow L^2(\Omega)$ is compact, and let  $V\in L^2_{\mathrm{loc}}(\Omega)$ be a potential taking only positive real values and satisfying the asymptotics $V(\cdot)\asymp\left\langle\cdot\right\rangle^{-\alpha}$ for some $\alpha\in\left]3,\infty\right[$.

For sufficient conditions on $\Omega$ guaranteeing the compact embedding, see the chapter 6 of \cite{AdamsEtFournier}, in particular Theorems 6.16 and 6.19, or the original article \cite{Adams}. The conditions are somewhat technical and therefore we do not reproduce them here. However, when $n\leqslant3$, one has the pleasant characterization: the embedding $H_0^2(\Omega)\longrightarrow L^2(\Omega)$ is compact if and only if $\Omega$ does not contain infinitely many pairwise disjoint balls which are all of the same size (see remarks 6.17.3, 6.9 and 6.11 in \cite{AdamsEtFournier}).

The theorems below cover as a special case bounded domains $\Omega$ with potentials $V$, which take only positive real values, and which are bounded and bounded away from zero.

In our setting, transmission eigenvalues for the Schr\"odinger equation are defined to be those complex numbers $\lambda$ for which there exist functions
\[v,w\in\bigl\{u\in H_{\mathrm{loc}}^2(\Omega)\bigm|\widetilde u\in B^\ast\bigr\}\setminus0\]
solving the equations
\[(-\Delta+V-\lambda)\,v=0,\quad\quad(-\Delta-\lambda)\,w=0\]
in $\Omega$,
and connected by the asymptotic relation and boundary conditions 
\[v-w\in\mathring B_2^\ast(\Omega)=\bigl\{u\in H_{\mathrm{loc}}^2(\Omega)\bigm|\widetilde u\in B_2^\ast(\mathbb R^n)\bigr\},\]
where $\widetilde u\colon\mathbb R^n\longrightarrow\mathbb C$ coincides with $u$ in $\Omega$ and vanishes identically elsewhere. It does no harm to occasionally identify $u$ with its zero extension $\widetilde u$.

The multiplicity of a transmission eigenvalue is defined as the dimension of the vector space of pairs of functions $\langle v,w\rangle$ solving the above problem.

We shall only consider real transmission eigenvalues and this is a genuine restriction (as was first shown by F.~Cakoni, D.~Colton and D.~Gintides \cite{CakoniEtColtonEtGintides}).

Our main theorem is

\begin{theorem}\label{main-theorem-1}
The set of positive real transmission eigenvalues for the Schr\"o\-ding\-er equation is a discrete subset of $\left[0,\infty\right[$, and each of its elements is of finite multiplicity.
\end{theorem}

For the Helmholtz equation the perturbed equation for $v$ is
\[(-\Delta+\lambda V-\lambda)\,v=0,\]
and one excludes the uninteresting value $\lambda=0$, but otherwise everything else is the same. In particular, we have

\begin{theorem}\label{main-theorem-2}
The set of positive real transmission eigenvalues for the Helmholtz equation is a discrete subset of $\left[0,\infty\right[$, and each of its elements is of finite multiplicity.
\end{theorem}

\section{Proof of Theorem~\ref{main-theorem-1}}

\subsection{Reduction to a fourth-order equation}

The first step in the proof is writing the transmission eigenvalue problem as a single fourth-order partial differential equation. This idea is rather standard and is the basis for many discreteness and existence proofs in the literature. The non-vanishing of $V$ is rather essential here.

We shall handle the operator in the fourth-order equation using quadratic forms, and this will require a shift from the $B^\ast$-based spaces to certain weighted $L^2$-based spaces.
The role of the ambient space will be played by $L_V$, a space which we define to consist of those $L^2_{\mathrm{loc}}$-functions in $\Omega$ whose zero extensions belong to Agmon's weighted space
\[L^{2,\alpha/2}(\mathbb R^n)=\left\{u\in L_{\mathrm{loc}}^2(\mathbb R^n)\,\middle|\vphantom{\Big|}\,\langle\cdot\rangle^{\alpha/2}u\in L^2(\mathbb R^n)\right\},\]
which of course is a Hilbert space with the right weighted $L^2$-norm.

The quadratic form domain will be $H_V$, a space which we define to consist of those $L_{\mathrm{loc}}^2$-functions in $\Omega$ whose zero extensions belong to Agmon's weighted space
\[H_{2,\alpha/2}(\mathbb R^n)=\left\{u\in L_{\mathrm{loc}}^2(\mathbb R^n)\,\vphantom{\Big|}\middle|\,
\partial^\gamma u\in L^{2,\alpha/2}(\mathbb R^n),\forall\left|\gamma\right|\leqslant2\right\}.\]
The space $H_V$ is Hilbert when equipped with the restriction of the $H_{2,\alpha/2}$-norm. 

We point out that $H_V$ embeds compactly into $L_V$. It is easy to split this embedding into three parts
\[H_V\longrightarrow H_0^2(\Omega)\longrightarrow L^2(\Omega)\longrightarrow L_V,\]
where the middle one is the obvious embedding, which is assumed to be compact, and the first and the last mappings are multiplications by $\langle\cdot\rangle^{\alpha/2}$ and $\langle\cdot\rangle^{-\alpha/2}$, respectively.

Now we are ready to state the transition to a fourth-order equation:

\begin{lemma}\label{Lemma}
If a positive real number $\lambda$ is a transmission eigenvalue then there exists a function $u\in H_V\setminus0$ solving the equation
\begin{equation}\label{fourth-order-equation}
(-\Delta+V-\lambda)\,\frac1V\,(-\Delta-\lambda)\,u=0
\end{equation}
in $\Omega$
in the sense of distributions.
Furthermore, this transition retains multiplicities.
\end{lemma}

If $v$ and $w$ solve the transmission eigenvalue problem, then it is a matter of simple calculation to see that $u=v-w$ solves the fourth-order equation.
It only remains to see that $u\in\mathring B_2^\ast(\Omega)$ corresponding to a transmission eigenvalue necessarily belongs to $H_V$. This follows from the observation that
\[\left(-\Delta-\lambda\right)u=-Vv,\]
not only in $\Omega$ but also in $\mathbb R^n$.
The function $\langle\cdot\rangle^{\alpha-1-\varepsilon}\,V$ has enough decay to be a short-range potential, and so
$\langle\cdot\rangle^{\alpha-1-\varepsilon}\,Vv\in B$. Now, by a basic inequality in short-range scattering theory (see e.g. Theorem 14.3.7 in \cite{HormanderII}),
\[\langle\cdot\rangle^{\alpha-1-\varepsilon}\,\partial^\gamma u\in B^\ast,\]
for $\left|\gamma\right|\leqslant2$. It is then easy to check that
\[\langle\cdot\rangle^{\alpha-3/2-\varepsilon}\,\partial^\gamma u\in L^2(\mathbb R^n),\]
again for $\left|\gamma\right|\leqslant2$, which in turn implies
\[\langle\cdot\rangle^{\alpha/2}\,\partial^\gamma u\in L^2(\mathbb R^n),\]
since $\alpha>3$.

\medbreak

From now on, we focus on studying the spectral properties of the fourth-order equation (\ref{fourth-order-equation}). In particular, we shall establish a discreteness result and a conditional existence result. The discreteness result, together with Lemma \ref{Lemma}, implies Theorem~\ref{main-theorem-1}.
The hypothesis required for the general existence result concerns the existence for suitable simple cases.
\begin{hypothesis}\label{the-hypothesis}
For any ball $B$ in $\mathbb R^n$ and any constant potential $V_0\in\mathbb R_+$, there exists a Schr\"odinger transmission eigenvalue.
\end{hypothesis}
\noindent We do not know whether this hypothesis is true.

By inspecting the conditional existence proof (which will be given in Section~\ref{existence-section}), we see that there is some freedom in the formulation of the hypothesis. For example, we only need the existence for a sequence of balls and positive constant potentials, where both the radii of the balls and the potentials tend to zero. Furthermore, balls could be replaced by any domains whose diameters shrink to zero, and the potentials do not have to be constant, as long as their $L^\infty$-norms tend to zero, and they are positive and bounded away from zero.

Also, it should be noted, that by considering radial functions (see Section~\ref{helmholtz-section}) one can prove that there are transmission eigenvalues for any ball, provided that $V_0$ is sufficiently large. From this the approach of Section~\ref{existence-section} will give unconditional existence of eigenvalues (not necessarily infinitely many), provided that the potential $V$ is sufficiently large in some balls in~$\Omega$.

\begin{theorem}
The set of real numbers $\lambda$ for which the equation (\ref{fourth-order-equation}) has a non-trivial $H_V$-solution is a discrete subset of $\left[0,\infty\right[$. For each such $\lambda$ the space of solutions is finite dimensional. Furthermore, if the Hypothesis~\ref{the-hypothesis} holds, the set of such real numbers $\lambda$ is infinite.
\end{theorem}

For a bounded $\Omega$ the $H_V$-solutions that can be conditionnally obtained by this theorem give rise to transmission eigenvalues, respecting multiplicities, and we get the existence of infinitely many transmission eigenvalues. Unfortunately, in the unbounded case this does not work. The obstacle is that the solutions to the fourth-order equation belong to weighted spaces which essentially guarantee that division by $\sqrt V$ is a reasonably good operation, whereas in order to get from the fourth-order equation back to the interior transmission problem one needs to divide by $V$, an operation genuinely worse than division by $\sqrt V$, and there seems to be no way of guaranteeing that the asymptotic behaviour of the apparent transmission eigenfunction pair is sufficiently good.

\subsection{The quadratic forms}

We will handle the operator on the left-hand side of (\ref{fourth-order-equation}) via quadratic forms, and for this purpose we define for each $\lambda\in\mathbb C$ the quadratic
form
\[Q_\lambda=u\longmapsto\left\langle(-\Delta+V-\overline\lambda)\,u\,\middle|\,\frac1V\,
(-\Delta-\lambda)\,u\right\rangle\colon H_V\longrightarrow\mathbb C,\]
where the $L^2$-inner product is linear in the second argument.
Instead of considering $Q_\lambda$ as a quadratic form in $L^2(\Omega)$, we shall consider it in the weighted $L^2$-space~$L_V$.
The idea of using weighted $L^2$-spaces as the ambient Hilbert spaces, in order to handle degenerate or even singular potentials in the case of a bounded domain, has been used in the papers \cite{Hickmann}, \cite{SerovEtSylvester} and \cite{Serov}, where the weight is a power of distance to the boundary of the domain.

It turns out that the family $\langle
Q_\lambda\rangle_{\lambda\in\mathbb C}$ has the pleasant
properties enumerated in the theorem below. An excellent reference for the basic theory of quadratic forms and analytic perturbation theory used is the book by Kato \cite{Kato}, in particular its Chapters VI and VII. More detailed references will be given in the course of the proofs of the statements.

\begin{theorem}
\begin{enumerate}
\item
The quadratic forms $Q_\lambda$ form an entire self-adjoint analytic family of forms of type (a) with compact resolvent, and therefore gives rise to a family of operators $T_\lambda$, which is an entire self-adjoint analytic family of operators of type (B) with compact resolvent.

\item
Furthermore, there exists a sequence $\left\langle\mu_\nu(\cdot)\right\rangle_{\nu=1}^\infty$ of real-analytic functions $\mu_\nu(\cdot)\colon\mathbb R\longrightarrow\mathbb R$ such that, for real $\lambda$, the spectrum of $T_\lambda$, which consists of a discrete set of real eigenvalues of finite multiplicity, consists of $\mu_1(\lambda)$, $\mu_2(\lambda)$, \dots, including multiplicity.

\item
In addition, for any given $T\in\mathbb R_+$, there exists constant $c\in\mathbb R_+$ such that
\[\bigl|\mu_\nu(\lambda)-\mu_\nu(0)\bigr|\ll_Te^{c|\lambda|}-1\]
for all $\lambda\in\left[-T,T\right]$ and each $\nu\in\mathbb Z_+$.

\item
The pairs $\left\langle\lambda,u\right\rangle\in\mathbb R\times H_V$ for which (\ref{fourth-order-equation}) holds,
are in bijective correspondence with the pairs $\left\langle\nu,\lambda\right\rangle\in\mathbb Z_+\times\mathbb R$ for which $\mu_\nu(\lambda)=0$.

\item If Hypothesis \ref{the-hypothesis} holds, then there are infinitely many such pairs $\left\langle\lambda,u\right\rangle$.
\end{enumerate}
\end{theorem}

The discreteness result follows easily from these properties of $Q_\lambda$. It is obvious
that zero is not an eigenvalue of $T_\lambda$ for any negative real $\lambda$,
as $Q_\lambda(u)>0$ for all non-zero functions $u\in\mathrm{Dom}\,Q_\lambda$. Hence none
of the functions $\mu_\nu(\cdot)$ can vanish identically, so that the set of
zeroes of each of them is discrete. Why the union of the zero sets can not have
an accumulation point follows immediately
from the third statement above,
which says that the functions $\mu_\nu(\cdot)$ change their values uniformly
locally exponentially. That is, when the value of $\lambda$ changes by a finite amount, only finitely many $\mu_\nu(\cdot)$ will have enough time to drop to zero.

The second statement follows immediately from a basic result in the perturbation
theory of linear operators, once the first has been proven; for this see
\cite[rem. VII.4.22, p.~408]{Kato} and the backwards references.
The third statement comes from theorem VII.4.21 \cite[p. 408]{Kato}.
The fourth and fifth statements will be consequences of the mini-max principle, but will
be given only after the first one has been dealt with.

We remark that the proof of the fifth statement only requires continuity of the family $\left\langle\mu_\nu(\cdot)\right\rangle_{\nu=1}^\infty$, which can be proved using the mini-max principle with no reference to non-real values of~$\lambda$ (see e.g. the proof of Lemma~12 in \cite{PaivarintaEtSylvester}). The observation that these eigenvalues depend real-analytically on~$\lambda$ seems to be new.

\subsection{A weighted inequality}

The proof of closedness of $Q_\lambda$ will depend on the following weighted inequality.

\begin{lemma}\label{L-weighted-inequality}
Let $K\subseteq\mathbb C$ be compact, and let $s\in\mathbb R$. Then
\[\bigl\|\langle\cdot\rangle^s\,u\bigr\|+\bigl\|\langle\cdot\rangle^s\,\nabla u\bigr\|+\bigl\|\langle\cdot\rangle^s\,\nabla\otimes\nabla u\bigr\|
\ll_{K,s}\bigl\|\langle\cdot\rangle^s\left(-\Delta-\lambda\right)u\bigr\|+\bigl\|\langle\cdot\rangle^s\,u\bigr\|\]
for all $u\in C_{\mathrm c}^\infty(\mathbb R^n)$ and $\lambda\in K$.
\end{lemma}

Here and elsewhere, given an expression $E(\cdot)$, we use the short-hand notations $E(\nabla)$ and $E(\nabla\otimes\nabla)$ for 
$\sum_{\left|\alpha\right|=1}E(\partial^\alpha)$
and
$\sum_{\left|\alpha\right|=2}E(\partial^\alpha)$,
respectively.
When necessary, we shall use other similar short-hands whose meaning will be clear.

\textit{Proof of Lemma~\ref{L-weighted-inequality}.} The following argument is an adaptation of the proof of Lemma A.3 of \cite[p.~206]{Agmon}. Since
\[\langle\cdot\rangle^4\ll_K\left|4\pi^2\left|\cdot\right|^2-\lambda\right|^2+1,\]
multiplication by $\bigl|\widehat u\bigr|^2$ and integration over $\mathbb R^n$ gives
\[\bigl\|u\bigr\|+\bigl\|\nabla u\bigr\|+\bigl\|\nabla\otimes\nabla u\bigr\|
\ll_K\bigl\|\left(-\Delta-\lambda\right)u\bigr\|+\bigl\|u\bigr\|.\]

In order to introduce weights, we observe that for $\varepsilon\in\left]0,1\right]$,
\[\langle\cdot\rangle\asymp_{\varepsilon,s}\langle\varepsilon\cdot\rangle,\]
and that
\[\partial^\alpha\langle\varepsilon\cdot\rangle^s\ll_{\varepsilon,s}\langle\varepsilon\cdot\rangle^s.\]
Now Leibniz's rule, the weightless inequality and the triangle inequality give
\begin{align*}
&\hspace{-1em}\bigl\|\langle\varepsilon\cdot\rangle^s\,u\bigr\|+\bigl\|\langle\varepsilon\cdot\rangle^s\,\nabla u\bigr\|+\bigl\|\langle\varepsilon\cdot\rangle^s\,\nabla\otimes\nabla u\bigr\|\\[1ex]
&\ll\bigl\|\langle\varepsilon\cdot\rangle^s\,u\bigr\|+\bigl\|\nabla(\langle\varepsilon\cdot\rangle^s\,u)\bigr\|+\bigl\|\nabla\otimes\nabla(\langle\varepsilon\cdot\rangle^s\,u)\bigr\|\\
&\qquad+\bigl\|(\nabla\langle\varepsilon\cdot\rangle^s)\,u\bigr\|+\bigl\|(\nabla\langle\varepsilon\cdot\rangle^s)\otimes\nabla u\bigr\|+\bigl\|(\nabla\otimes\nabla\langle\varepsilon\cdot\rangle^s)\,u\bigr\|\\[1ex]
&\ll_K\bigl\|(-\Delta-\lambda)\,(\langle\varepsilon\cdot\rangle^s\,u)\bigr\|+\bigl\|\langle\varepsilon\cdot\rangle^s\,u\bigr\|\\
&\qquad+\bigl\|(\nabla\langle\varepsilon\cdot\rangle^s)\,u\bigr\|+\bigl\|(\nabla\langle\varepsilon\cdot\rangle^s)\otimes\nabla u\bigr\|+\bigl\|(\nabla\otimes\nabla\langle\varepsilon\cdot\rangle^s)\,u\bigr\|\\[1ex]
&\ll\bigl\|\langle\varepsilon\cdot\rangle^s\,(-\Delta-\lambda)\,u\bigr\|
+\bigl\|(\nabla\langle\varepsilon\cdot\rangle^s)\cdot\nabla u\bigr\|+\bigl\|(\Delta\langle\varepsilon\cdot\rangle^s)\,u\bigr\|+\bigl\|\langle\varepsilon\cdot\rangle^s\,u\bigr\|\\
&\qquad+\bigl\|(\nabla\langle\varepsilon\cdot\rangle^s)\,u\bigr\|+\bigl\|(\nabla\langle\varepsilon\cdot\rangle^s)\otimes\nabla u\bigr\|+\bigl\|(\nabla\otimes\nabla\langle\varepsilon\cdot\rangle^s)\,u\bigr\|\\[1ex]
&\ll_s\bigl\|\langle\varepsilon\cdot\rangle^s\,(-\Delta-\lambda)\,u\bigr\|
+\varepsilon\bigl\|\langle\varepsilon\cdot\rangle^s\,\nabla u\bigr\|+\varepsilon^2\bigl\|\langle\varepsilon\cdot\rangle^s\,u\bigr\|+\bigl\|\langle\varepsilon\cdot\rangle^s\,u\bigr\|\\
&\qquad+\varepsilon\bigl\|\langle\varepsilon\cdot\rangle^s\,u\bigr\|+\varepsilon\bigl\|\langle\varepsilon\cdot\rangle^s\,\nabla u\bigr\|+\varepsilon^2\bigl\|\langle\varepsilon\cdot\rangle^s\,u\bigr\|\\[1ex]
&\ll\bigl\|\langle\varepsilon\cdot\rangle^s\,(-\Delta-\lambda)\,u\bigr\|+\bigl\|\langle\varepsilon\cdot\rangle^s\,u\bigr\|+\varepsilon\bigl\|\langle\varepsilon\cdot\rangle^s\,\nabla u\bigr\|.
\end{align*}
Choosing a sufficiently small $\varepsilon$, subject to the choices of $K$ and $s$, allows us to eliminate the first-order term from the right-hand side, giving the weighted version of the desired inequality.

\subsection{$Q_\lambda$ is a good self-adjoint family}

The fact that the family $\langle T_\lambda\rangle_{\lambda\in\mathbb C}$
is a self-adjoint analytic family of type (B) with compact resolvent will
follow from a number of different results in the aforementioned book
\cite{Kato}.

If $\langle Q_\lambda\rangle$ form a self-adjoint analytic family of
quadratic forms of type (a), then for each $\lambda\in\mathbb C$, there
corresponds a unique closed linear operator $T_\lambda$; since $Q_\lambda$ is
densely defined, sectorial and closed (as will be shown later), the unique
existence of $T_\lambda$ is given by \cite[thm. VI.2.1, p. 322]{Kato}, and
the operator $T_\lambda$ is furthermore $m$-sectorial.

The theorem VII.4.2 \cite[p. 395]{Kato} then says that $\langle
T_\lambda\rangle$ is an analytic family of operators (in the sense of Kato).
Since
\[\mathrm{Dom}\,T_\lambda\subseteq\mathrm{Dom}\,Q_\lambda=H_V,\]
and $H_V$ embeds compactly into $L_V$, the family $\langle
T_\lambda\rangle$ has compact resolvent. Finally, the family is self-adjoint,
i.e. $T_\lambda^\ast=T_{\overline\lambda}$, since $\overline{Q_\lambda}
=Q_{\overline\lambda}$ for all $\lambda\in\mathbb C$.
This follows from theorem VI.2.5 \cite[p. 323]{Kato}.
In particular, $T_\lambda$ is a self-adjoint operator with compact resolvent
for real $\lambda$.

Thus it remains to prove that $\langle Q_\lambda\rangle$ is an analytic family
of type (a). By definition, this entails checking that
\begin{itemize}\setlength{\itemsep}{0pt}
\item Each $Q_\lambda$ is sectorial and closed, and $\mathrm{Dom}\,Q_\lambda$
is independent of $\lambda$; and
\item $Q_\lambda(u)$ is an entire function of $\lambda$ for any fixed $u\in
\mathrm{Dom}\,Q_\lambda$.
\end{itemize}

The latter condition is obviously satisfied as $Q_\lambda(u)$ is, in fact,
a second degree polynomial in $\lambda$. That $\mathrm{Dom}\,Q_\lambda$ is independent
of $\lambda$ is also obvious here, because the domain is simply $H_V$.
So it only remains to prove that each $Q_\lambda$ is sectorial and closed.

That $Q_\lambda$ is sectorial simply means that the set
$Q_\lambda[\{u\in H_V\,|\,\|u\|_{L_V}=1\}]$ is contained in a
sector-shaped set of the form
\[\{z\in\mathbb C\,|\,\arg(z-z_0)\leqslant\vartheta\}\]
for some fixed $z_0\in\mathbb C$ and $\vartheta\in\left[0,\frac\pi2\right[$.
This sectoriality condition can be established by the usual elementary arguments; see e.g. Example 1.7 in \cite[p.~312]{Kato}.

That $Q_\lambda$ is closed follows now from the fact that, by the weighted inequality proved above, the $H_{2,\alpha/2}$-norm and the norm arising from $Q_\lambda$, given by the expression
\[\sqrt{\Re Q_\lambda(\cdot)+(1+\lambda)\|\cdot\|_{L_V}^2},\]
are comparable on $C_{\mathrm c}^\infty(\Omega)$ and therefore the domain of $Q_\lambda$ is really just the closure of test functions of $\Omega$ in the right norm.

\subsection{The bijective correspondence between $\left\langle\lambda,u\right\rangle$ and $\left\langle\nu,\lambda\right\rangle$}

If zero is an eigenvalue of $T_\lambda$ with an eigenfunction $u\in H_V$,
then clearly $Q_\lambda(v,u)=0$ for all $v\in C_{\mathrm c}^\infty(\Omega)$, and $u$ is a non-trivial solution to the equation (\ref{fourth-order-equation}).

The other direction is only slightly more challenging to establish. Suppose that
for $\lambda\in\mathbb R$ the equation (\ref{fourth-order-equation}) has a non-trivial space of solutions in $H_V$ of dimension $N$.
Then $Q_\lambda$ vanishes in some subspace $Y\subseteq H_V$ of
dimension $N$, and in fact,
$Q_\lambda(u,v)=0$ for all $v\in H_V$ and $u\in Y$.
Our goal is to prove that zero is an
eigenvalue of $T_\lambda$ of multiplicity at least $N$
using the mini-max principle. Let the spectrum of $T_\lambda$ be
\[\mu_1\leqslant\mu_2\leqslant\mu_3\leqslant\dots.\]

The space $X$ corresponding to the negative eigenvalues of $T_\lambda$ is
finite dimensional, say of dimension $m\geqslant0$. Now the restriction
\[T|_{X^\perp}\colon X^\perp\cap\mathrm{Dom}\,T_\lambda\longrightarrow
X^{\perp}\]
is again a self-adjoint operator with compact resolvent and no negative
eigenvalues. The eigenvalues $\mu_{m+1}$, $\mu_{m+2}$, \dots, $\mu_{m+N}$
all have to be non-negative.

Conversely, $\mu_{m+N}$ is at most
\begin{align*}
&\max
\left\{Q_\lambda(f)\,\middle|\,\vphantom{\Big|}f\in\mathrm{span}\,\{X,Y\},\|f\|_{L_V}=1\right\}\\
=&\max
\left\{\left(Q_\lambda(g)
+2\Re Q_\lambda(g,h)+Q_\lambda(h)\right)\,\middle|\,\vphantom{\Big|}g\in X,h\in Y,\|g+h\|_{L_V}=1\right\},
\end{align*}
and in the expression $(\ldots)$ the first term is certainly $\leqslant0$ and
the remaining terms vanish. Thus
$\mu_{m+1}\leqslant\mu_{m+2}\leqslant\ldots\leqslant\mu_{m+N}\leqslant0$ and we are done.

\subsection{The conditional infinitude of zeroes of $\mu_\nu(\cdot)$}\label{existence-section}

Next we shall prove that, under Hypothesis~\ref{the-hypothesis}, for arbitrarily large positive integers~$N$, there exists at least $N$ pairs $\left\langle\nu,\lambda\right\rangle\in\mathbb Z_+\times\mathbb R$ satisfying $\mu_\nu(\lambda)=0$. This is achieved by comparison to the simpler domains with constant potentials for which the existence of a single transmission eigenvalue is guaranteed by Hypothesis~\ref{the-hypothesis}.

We choose $N$ small balls $B_1$, $B_2$, \dots, $B_N$, whose closures are in $\Omega$ and pairwise disjoint, and consider on them a constant potential $V_0\in\mathbb R_+$ such that $V_0\leqslant V$ in $B_1\cup B_2\cup\ldots\cup B_N$, and such that there is a number $\lambda\in\mathbb R_+$ which is a transmission eigenvalue for each of the balls. The above theorem guarantees the existence of such a small $V_0$.

The $H_0^2$-spaces of the balls naturally embed into $H_V$ by taking zero extensions of their elements. Denote by~$H(N)$ the closed subspace spanned by the images of the differences of the transmission eigenfunction pairs of $V_0$ in the small balls. This space has dimension at least~$N$.

Now the quadratic form $\widetilde Q_\lambda$ corresponding to the constant potential $V_0$ in $\Omega$ is basically the $Q_\lambda$ defined above, but with $1/V$ replaced by $1/V_0$. In particular, we have the inequality $Q_\lambda\leqslant\widetilde Q_\lambda=0$ in $H(N)$. (The domain of $\widetilde Q_\lambda$ can be chosen to be anything reasonable that contains $H(N)$ as we only need this non-positivity inequality.)

The eigenvalues of $T_\kappa$ are positive for $\kappa\in\mathbb R_-$, but by the mini-max principle, at least $N$ of the eigenvalues of $T_\lambda$ are non-positive. Therefore the functions $\mu_n(\cdot)$ must have at least $N$ zeroes in the interval $\left[0,\lambda\right]$.

\section{Some remarks on the Helmholtz case}\label{helmholtz-section}

Everything we do works for the Helmholtz equation with modest modifications. The fourth-order equation (\ref{fourth-order-equation}) should be replaced by
\begin{equation}\label{fourth-order-helmholtz-equation}
\left(-\Delta+\lambda V-\lambda\right)\frac1V\left(-\Delta-\lambda\right)u=0,
\end{equation}
and the quadratic forms should be redefined accordingly. The spectral properties will in this case be slightly better than in the Schr\"odinger case:
\begin{theorem}
The set of positive real numbers $\lambda$ for which the equation~(\ref{fourth-order-helmholtz-equation}) has a non-trivial $H_V$-solution is an infinite discrete subset of $\left[0,\infty\right[$, and for each such $\lambda$ the space of solutions is finite dimensional. Furthermore, the number of such $\lambda$ not exceeding $x\in\mathbb R_+$, counting multiplicities, is $\gg x^{n/2}$ as $x\longrightarrow\infty$.
\end{theorem}

The unconditional existence proof depends on
\begin{theorem}\label{concrete-existence}
For a ball $B$ in $\mathbb R^n$, and an arbitrarily small constant potential $c\in\mathbb R_+$, there exist infinitely many positive real Helmholtz transmission eigenvalues.
\end{theorem}
\noindent This is a special case of a much more general theorem on existence for radial potentials, a proof of which in three dimensions may be found in
\cite[p.~16]{ColtonEtPaivarintaEtSylvester}. For constant potentials, the proof simplifies nicely, and though it seems that there is no $n$-dimensional proof in the literature, the $3$-dimensional proof generalizes easily: the crucial difference is that $j_0(r)$ must be replaced by $r^{(2-n)/2}\,J_{(n-2)/2}(r)$.

Now we do not immediately see that $\mu_n(\lambda)>0$ for negative $\lambda$ and each $n\in\mathbb Z_+$. Instead, we observe easily that $\mu_n(0)>0$ for each $n$: If $Q_0(u)=0$, then $\Delta u\equiv0$, implying that $\nabla u\equiv0$, and therefore $u$ must vanish.

The asymptotic lower bound $\gg x^{n/2}$ for the number of
zeroes of $\mu_\nu(\cdot)$
not exceeding a large positive real number~$x$ follows from the fact that transmission eigenvalues for the Helmholtz equation scale under dilations like Dirichlet eigenvalues.

More precisely, let us look at a ball $B$ whose closure is contained in $\Omega$, and let $V_0\in\mathbb R_+$ be so small that $V_0\leqslant V(\cdot)$ in $B$. Now there exists a transmission eigenvalue $\lambda$ for $B$ and the constant potential $V_0$.

Given any $\varepsilon\in\mathbb R_+$, it is easy to see that the number $\frac\lambda{\varepsilon^2}$ is a transmission eigenvalue for any translate of $\varepsilon B$ with the constant potential $V_0$.

Let $x\in\mathbb R_+$, and choose $\varepsilon=\lambda^{1/2}\,x^{-1/2}$. Now the number $\frac\lambda{\varepsilon^2}=x$ is a transmission eigenvalue for any translate of $\varepsilon B$ with the constant potential $V_0$, and we can pack $\gg\varepsilon^{-n}\gg_\lambda x^{n/2}$ such translates inside $B$ so that no two of them intersect, provided that $x$ is large enough. These will correspond to the balls $B_1$, $B_2$, \dots, $B_N$ of section~\ref{existence-section}. Now we finish the proof in the same way as in section~\ref{existence-section} and obtain $\gg_\lambda x^{n/2}$ zeroes not exceeding~$x$.

\section*{Acknowledgements}

This research was funded by Finland's Ministry of Education through the Doctoral Program in Inverse Problems, and by the Finnish Centre of Excellence in Inverse Problems Research.
The author wishes to thank Prof. M. Salo for suggesting this research topic and for a number of helpful and encouraging discussions.

\end{document}